\documentclass[a4paper,11pt, reqno, oneside]{amsart}

\usepackage{amssymb}
\usepackage{amsmath}
\usepackage{amstext}
\usepackage{amscd}
\usepackage{amsthm}
\usepackage{amsfonts}
\usepackage{enumerate}
\usepackage{latexsym}
\usepackage{mathrsfs}
\usepackage{psfrag}
\usepackage{epsfig}
\usepackage{epic}
\usepackage{mathtools}
\usepackage{array}
\usepackage{color}
\usepackage[all]{xy}
\usepackage{bm}
\usepackage{caption}

\makeatletter
  
  \@addtoreset{equation}{section}
\makeatother

\newtheorem{theorem}{Theorem}[section]
\newtheorem{lemma}[theorem]{Lemma}
\newtheorem{proposition}[theorem]{Proposition}
\newtheorem{corollary}[theorem]{Corollary}

\theoremstyle{definition}
\newtheorem{definition}[theorem]{Definition}
\newtheorem{definiton-theorem}[theorem]{Definition-Theorem}
\newtheorem{lemma-definition}[theorem]{Lemma-Definition}
\newtheorem{example}[theorem]{Example}

\theoremstyle{remark}
\newtheorem{remark}[theorem]{Remark}

\newcommand{\calF}{\mathcal{F}}
\newcommand{\calG}{\mathcal{G}}

\newcommand{\calN}{\mathcal{N}}

\newcommand{\fkm}{\mathfrak{m}}

\newcommand{\fkS}{\mathfrak{S}}

\newcommand{\rad}{\operatorname{rad}}
\newcommand{\rank}{\operatorname{rank}}
\newcommand{\Ext}{\operatorname{Ext}}

\newcommand{\tor}{\operatorname{tor}}
\newcommand{\Hom}{\operatorname{Hom}}

\newcommand{\Ker}{\operatorname{Ker}}

\newcommand{\GL}{\operatorname{GL}}
\newcommand{\SL}{\operatorname{SL}}

\newcommand{\charac}{\operatorname{char}}

\newcommand{\Irr}{\operatorname{Irr}}

\newcommand{\twoheadlongrightarrow}{\relbar\joinrel\twoheadrightarrow}

\begin{document}

%%%%%---title---%%%%%
\title[Dual $F$-signature of special CM modules]{Dual $F$-signature of special Cohen-Macaulay modules over cyclic quotient surface singularities}
\author[Yusuke~Nakajima]{YUSUKE NAKAJIMA}
\date{}


\address{Graduate School Of Mathematics, Nagoya University, Furocho, Chikusaku, Nagoya, 464-8602 Japan} 
\email{m06022z@math.nagoya-u.ac.jp}
\maketitle

%%%---abstract--------------------------------------------------------------------------------------------
\begin{abstract}
The notion of $F$-signature was defined by C.~Huneke and G.~Leuschke and this numerical invariant characterizes some singularities. 
This notion is extended to finitely generated modules by A.~Sannai and is called dual $F$-signature.  
In this paper, we determine the dual $F$-signature of a certain class of Cohen-Macaulay modules (so-called ``special") over cyclic quotient surface singularities. 
Also, we compare the dual $F$-signature of a special Cohen-Macaulay module with that of its Auslander-Reiten translation. 
This gives a new characterization of the Gorensteiness. 
\end{abstract}

%\tableofcontents

%%%---text start------------------------------------------------------------------------------------------
\section{Introduction}
\label{intro}

Throughout this paper, we suppose that $k$ is an algebraically closed field of prime characteristic $p>0$.
Let $(R,\fkm,k)$ be a Noetherian local ring with $\operatorname{char} R=p>0$. 
Since $\operatorname{char}R=p>0$, we can define the Frobenius map $F:R\rightarrow R\;(r\mapsto r^p)$. 
For $e\in\mathbb{N}$, we also define the $e$-times iterated Frobenius map $F^e:R\rightarrow R\;(r\mapsto r^{p^e})$. 
For any $R$-module $M$, we denote the module $M$ with its $R$-module structure pulled back via the $e$-times iterated Frobenius map $F^e$ by ${}^eM$.
Namely, ${}^eM$ is just $M$ as an abelian group, and its $R$-module structure is defined by $r\cdot m\coloneqq F^e(r)m=r^{p^e}m$ for all $r\in R,\;m\in M$. 
We say $R$ is $F$-finite if ${}^1R$ is a finitely generated $R$-module.

\medskip

In order to investigate the properties of $R$, C.~Huneke and G.~Leuschke introduced the notion of $F$-signature.

\begin{definition}[\cite{HL}]
\label{HLFsig}
Let $(R,\fkm,k)$ be a reduced $F$-finite local ring of prime characteristic $p>0$.
For each $e\in\mathbb{N}$, decompose ${}^eR$ as follows
\[
{}^eR\cong R^{\oplus a_e}\oplus M_e,
\]
where $M_e$ has no free direct summands. We call $a_e$ the $e$-th $F$-splitting number of $R$.
Then, the $F$-signature of $R$ is 
\[
 s(R)\coloneqq\lim_{e\rightarrow\infty}\frac{a_e}{p^{ed}},
\]
if it exists, where $d\coloneqq\dim R$.
\end{definition}

Note that K.~Tucker showed its existence in a general situation \cite{Tuc}. 
As Kunz's theorem \cite{Kun} shows, this invariant measures the deviation from regularity $($see also Theorem~\ref{charac_sing}$(1)\,)$. 

   \medskip
For a finitely generated $R$-module, A.~Sannai extended the notion of $F$-signature as follows.

\begin{definition}[\cite{San}]
\label{SanFsig}
Let $(R,\fkm,k)$ be a reduced $F$-finite local ring of prime characteristic $p>0$.
For a finitely generated $R$-module $M$ and $e\in\mathbb{N}$, we set
\[
 b_e(M)\coloneqq\operatorname{max}\{n\;|\;\exists\varphi:{}^eM\twoheadrightarrow M^{\oplus n}\}, 
\]
and call it the $e$-th $F$-surjective number of $M$. 
Then we call the limit
\[
 s(M)\coloneqq\lim_{e\rightarrow\infty}\frac{b_e(M)}{p^{ed}}
\]
the dual $F$-signature of $M$ if it exists, where $d=\dim R$.
\end{definition}

\begin{remark}
The dual $F$-signature of $R$ coincides with the $F$-signature of $R$, because the morphism 
${}^eR\twoheadrightarrow R^{\oplus b_e(R)}$ is split. Therefore, we use the same notation unless it causes confusion. 
\end{remark}

By using these invariants, we can characterize some singularities. 

\begin{theorem}[\cite{HL}, \cite{Yao2}, \cite{AL}, \cite{San}]
\label{charac_sing}
Let $(R,\fkm,k)$ be a $d$-dimensional reduced $F$-finite Noetherian local ring with $\charac R=p>0$. Then we obtain 
\begin{itemize}
 \item [(1)] $R$ is regular if and only if $s(R)=1$, 
 \item [(2)] $R$ is strongly $F$-regular if and only if $s(R)>0$.
\end{itemize}

In addition, we suppose $R$ is Cohen-Macaulay with the canonical module $\omega_R$, then 
\begin{itemize}
 \item [(3)] $R$ is $F$-rational if and only if $s(\omega_R)>0$, 
 \item [(4)] $s(R)\le s(\omega_R)$,
 \item [(5)] $s(R)= s(\omega_R)$ if and only if $R$ is Gorenstein.
\end{itemize}
\end{theorem}

As the above theorem shows, the value of $s(R)$ and $s(\omega_R)$ contain some information regarding singularities. 
How about the value of the dual $F$-signature for other $R$-modules ? 
The value of the dual $F$-signature is not known except in the case of two-dimensional Veronese subrings \cite{San}. 
%and we don't know the method for determining it, this question is so difficult for now. 
Therefore, in this paper, we determine the dual $F$-signature for a certain class of Cohen-Macaulay (= CM) modules (so-called special CM modules) 
over cyclic quotient surface singularities.  
%As we will see later, special CM modules are compatible with the geometry. 

The study of special CM modules was started by the work of J.~Wunram \cite{Wun1}, \cite{Wun2} (the definition of special CM modules appears in Section~\ref{section_DFsig}).
For a finite subgroup $G$ of $\SL(2,k)$ such that the order of $G$ is invertible in $k$, the McKay correspondence is very famous, that is, there is a one-to-one correspondence between 
non-trivial irreducible representations of $G$ and irreducible exceptional curves on the minimal resolution of quotient surface singularity. 
When we intend to generalize this correspondence to a finite subgroup $G$ of $\GL(2,k)$, this correspondence is no longer true. 
In fact, there are more irreducible representations than exceptional curves. However, if we choose some irreducible representations which are called special, 
then we again obtain one-to-one correspondence between irreducible special representations of $G$ and exceptional curves \cite{Wun2}, 
and a maximal CM module associated with a special representation is called a special CM module.
For more about the special McKay correspondence, see also \cite{Ish},\cite{Ito} and \cite{Rie}.

\begin{remark}
\label{remSL2}
All irreducible representations of a finite subgroup of $\SL(2,k)$ are special, 
thus we can recover the McKay correspondence in the original sense from the special one. 
\end{remark}

For a cyclic quotient singularity, a special CM module takes the following simple form. 
(For more details on terminologies, see Section~\ref{Preliminary} and \ref{section_DFsig}.)

Suppose $R$ is the invariant subring of $S=k[[x, y]]$ under the action of a cyclic group $\frac{1}{n}(1,a)$. 
In this situation, a non-free indecomposable special CM $R$-module is described as $M_{i_t}=Rx^{i_t}+Ry^{j_t}$ 
(i.e. it is minimally 2-generated). 
%Then we have the value of the dual $F$-signature as follows. Especially, they are all rational. 
The following theorem gives the value of the dual $F$-signature; note that they are all rational. 

\begin{theorem}[see Theorem~\ref{main}]
For any non-free indecomposable special CM $R$-module $M_{i_t}$, we have 
\begin{equation*}
s(M_{i_t})=\begin{cases} \displaystyle\frac{\operatorname{min}(i_t,j_t)+1}{n}&\text{(if \;$i_t\neq j_t$)} \\
                         &\\
                         \displaystyle\frac{2i_t+1}{2n}&\text{(if \;$i_t=j_t$)}. 
           \end{cases}
\end{equation*}

\end{theorem}

\bigskip

Moreover, by paying attention to special CM modules and their Auslander-Reiten translations, 
we characterize when the ring is Gorenstein. 

\begin{theorem}[see Theorem~\ref{comp_AR_tr}]
\label{intro_comp_AR_tr}
Let $R$ be a quotient surface singularity. 
(Remark that we don't restrict to a cyclic case.) 
Suppose $M$ is an indecomposable special CM $R$-module. Then we have 
\[
s(M)\le s(\tau(M)).
\]
Moreover, if $s(M)= s(\tau(M))$ for an indecomposable special CM $R$-module $M$, then $R$ is Gorenstein. 
(Note that if $R$ is Gorenstein, then $s(M)=s(\tau(M))$ holds for all indecomposable MCM modules.) 
\end{theorem}

\begin{remark}
Since $\tau(R)\cong\omega_R$ in our situation, this theorem is an analogue of 
Theorem~\ref{charac_sing}\,$(4), (5)$. 
But it says that this characterization is obtained by not only the comparison between $R$ and $\omega_R$ but also 
the comparison between a special CM module and its AR translation.
\end{remark}

The structure of this paper is as follows.
In order to determine the dual $F$-signature, we need the notion of the generalized $F$-signature and the Auslander-Reiten quiver. 
Thus, we prepare them in Section~\ref{Preliminary}. In Section~\ref{section_DFsig}, we determine the dual $F$-signature of special CM modules over cyclic quotient surface singularities and give several examples. 
In Section~\ref{compare_AR_translation}, we compare a special CM module with its Auslander-Reiten translation by using the dual $F$-signature 
and characterize the Gorensteiness. Note that the statements appearing in Section~\ref{compare_AR_translation} hold not only for 
cyclic quotient surface singularities but also for any quotient surface singularities.

%%%%%%%%%%%%%%%%%%%%%%%%%%%%%%%%%%%%%%%%%%%%%%%%%%%%%%%%%%%%%%%%%%%%%%%%%%%%%%%%%%%%%%%%%%%%%%%%%%%%%%%%%%%%
\section{Preliminary}
\label{Preliminary}
\subsection{Generalized $F$-signature of invariant subrings}
\label{subsec_genFsig}

Let $G$ be a finite subgroup of $\GL(d,k)$ which contains no pseudo-reflections except the identity, and $S\coloneqq k[[x_1,\cdots,x_d]]$ be a power series ring. 
We assume that the order of $G$ is coprime to $p=\mathrm{char}\;k$. We denote the invariant subring of $S$ under the action of $G$ by $R\coloneqq S^G$.
In order to determine the dual $F$-signature of a finitely generated $R$-module $M$, we have to know about the structure of ${}^eM$ 
(for instance, the direct sum decomposition of ${}^eM$, 
the asymptotic behavior of the multiplicities of direct summands, etc). To achieve this, we use the results of the generalized $F$-signature of invariant subrings \cite{HN}.

   \medskip
For a positive characteristic Noetherian ring, K.~Smith and M.~Van den Bergh introduced the notion of finite $F$-representation type \cite{SVdB}.
This notion is a characteristic $p$ analogue of the notion of finite representation type. 
The definition of finite $F$-representation type is the following.

\begin{definition}[\cite{SVdB}] 
We say that $R$ has finite $F$-representation type $($or FFRT for short$)$ by $\calN$ if 
there exists a finite set $\calN$ of isomorphism classes of indecomposable finitely generated $R$-modules,
such that for every $e\in\mathbb{N}$, the $R$-module ${}^eR$ is isomorphic to a finite direct sum of elements of $\calN$. 
\end{definition}

For example, a power series ring $S$ has FFRT by $\{S\}$ (cf. Kunz's theorem \cite{Kun}) and FFRT is inherited by a direct summand \cite{SVdB}. 
Thus, an invariant subring $R$ also has FFRT. More explicitly, we have the next proposition.

\begin{proposition}[\cite{SVdB}]
Let $V_0=k,V_1,\cdots,V_{n-1}$ be the complete set of irreducible representations of $G$ and 
we set $M_t\coloneqq(S\otimes_kV_t)^G\;\;(t=0,1,\cdots,n-1)$.   
Then $R$ has finite $F$-representation type by the finite set $\{M_0\cong R,M_1,\cdots,M_{n-1}\}$.
\end{proposition}

Thus we can write ${}^eR$ as follows.
\[
{}^eR\cong R^{\oplus c_{0,e}}\oplus M_1^{\oplus c_{1,e}}\oplus\cdots\oplus M_{n-1}^{\oplus c_{n-1,e}}.
\]

\begin{remark}
We can see that each $M_t$ is an indecomposable maximal Cohen-Macaulay (= MCM) $R$-module and $M_s\not\cong M_t\; (s\neq t)$ 
under the assumption $G$ contains no pseudo-reflections except the identity. 
Also, the multiplicities $c_{i,e}$ are determined uniquely in that case. For more details, we refer the reader to \cite[Section\,2]{HN}.
\end{remark}

Moreover, since an invariant subring $R$ has FFRT, the limit 
$\displaystyle\lim_{e\rightarrow\infty}\frac{c_{t,e}}{p^{de}}\;(t=0,1,\cdots,n-1)$ exists \cite{SVdB}, \cite{Yao}. 
Therefore we can define the limit $s(R,M_t)\coloneqq\displaystyle\lim_{e\rightarrow\infty}\frac{c_{t,e}}{p^{de}}$ and call it the generalized $F$-signature of $M_t$ 
with respect to $R$. 
The value of $s(R,M_t)$ is determined by M.~Hashimoto and the author as follows.

\begin{theorem}[\cite{HN}]
\label{gene-Fsig}
For $t=0,1,\cdots,n-1$, one has 
\[
s(R,M_t)=\frac{\rank_RM_t}{|G|}
\]
\end{theorem}

\begin{remark}
In the case of $t=0$ is also due to \cite{WY}.
A similar result holds for a finite subgroup scheme of $\mathrm{SL_2}$ \cite{HS}. 
\end{remark}

We also obtain the next statement as a corollary.

\begin{corollary}[\cite{HN}]
\label{gene-Fsig-cor}
Suppose an MCM $R$-module $M_t$ decomposes as
\[
{}^eM_t\cong R^{\oplus d^t_{0,e}}\oplus M_1^{\oplus d^t_{1,e}}\oplus\cdots\oplus M_{n-1}^{\oplus d^t_{n-1,e}}.
\]
Then, for all $t,u=0,\cdots,n-1$, we obtain 
\[
s(M_t,M_u)\coloneqq\lim_{e\rightarrow\infty}\frac{d^t_{u,e}}{p^{de}}=(\rank_RM_t)\cdot s(R,M_u)=\frac{(\rank_RM_t)\cdot(\rank_RM_u)}{|G|}.
\]
\end{corollary}

\begin{remark}
In dimension two, it is known that an invariant subring $R$ is of finite representation type, that is, it has only finitely many 
non-isomorphic indecomposable MCM $R$-modules $\{R,M_1,\cdots,M_{n-1} \}$. 
From Corollary~\ref{gene-Fsig-cor}, every indecomposable MCM $R$-module appears in ${}^eM_t$ as a direct summand for sufficiently large $e$. 
Thus, the additive closure $\mathrm{add}_R({}^eM_t)$ coincides with the category of MCM $R$-modules $\mathrm{CM}(R)$. 
So we apply several results in Auslander-Reiten theory to $\mathrm{add}_R({}^eM_t)$ (see the next subsection).
\end{remark}

%%%%%%%%%%%%%%%%%%%%%%%%%%%%%%%%%%%%%%%%%%%%%%%%%%%%%%%%%%%%%%%%%%%%%%%%%%%%%%%%%%%%%%%%%%%%%%%%%%%%%%%%%%%%
\subsection{Auslander-Reiten quiver}
\label{subsec_AR}

In this subsection, we review some results of Auslander-Reiten theory.
For details, see some textbooks (e.g.\cite{LW},\;\cite{Yo}) or \cite{Aus1}, \cite{Aus2}. 
We only discuss such a theory for the case of an invariant subring $R=S^G$ in $\dim R=2$. 

\begin{definition}[Auslander-Reiten sequence]
Let $R$ be an invariant subring and $M,N$ be indecomposable MCM $R$-modules.
We call a non-split short exact sequence
\[
0\rightarrow N\overset{f}{\rightarrow}L\overset{g}{\rightarrow}M\rightarrow 0
\]
the Auslander-Reiten $($= AR$)$ sequence ending in $M$ if for all MCM modules $X$ and for any morphism $\varphi:X\rightarrow M$ 
which is not a split surjection there exists $\phi:X\rightarrow L$ such that $\varphi=g\circ\phi$.
\end{definition}

Since $R$ is an isolated singularity, 
there exists the AR sequence ending in $M_t$ for each non-free indecomposable MCM $R$-module $M_t$, and it is unique up to isomorphism \cite{Aus2}. 
Concretely, the AR sequence ending in $M_t \;(t\neq 0)$ is 
\[
0\longrightarrow (S\otimes_k(\wedge^2V\otimes_kV_t))^G\longrightarrow (S\otimes_k(V\otimes_kV_t))^G\longrightarrow M_t=(S\otimes_kV_t)^G\longrightarrow 0, 
\]
where $V$ is a natural representation of $G$. 

In the case of $t=0$, there is the exact sequence 
\[
0\longrightarrow\omega_R=(S\otimes_k\wedge^2V)^G\longrightarrow(S\otimes_kV)^G\longrightarrow R=S^G\longrightarrow k\longrightarrow 0.
\]
This exact sequence is called the fundamental sequence of $R$.

We call the left term of these sequences the Auslander-Reiten (= AR) translation of $M_t$, and denote it by $\tau(M_t)$. 
Sometimes we denote the middle term of the AR sequence by $E_{M_t}$. 
It is known that $\tau(M_t)\cong(M_t\otimes_R\omega_R)^{**}$, where $(-)^*=\Hom_R(-,R)$ is the $R$-dual functor \cite{Aus1}. 
Note that $\tau(M_t)=M_{t-a-1}$ and $E_{M_t}=M_{t-1}\oplus M_{t-a}$ for $t=0,1,\cdots,n-1$ in the case of subsection~\ref{subsec_cyclic}. 

   \medskip
Next, we prepare some notions to define the Auslander-Reiten quiver. 

\begin{definition}[Irreducible morphism]
Let $M$ and $N$ be MCM $R$-modules. We decompose $M$ and $N$ into indecomposable modules as 
$M=\oplus_iM_i$,\;$N=\oplus_jN_j$ and decompose $\psi\in \Hom_R(M,N)$ along this decomposition as 
$\psi=(\psi_{ij}:M_i\rightarrow N_j)_{ij}$. 
Then we define the submodule $\rad_R(M,N)\subset\Hom_R(M,N)$ as follows .
\[
 \psi\in\rad_R(M,N)\overset{def}{\Longleftrightarrow} \text{no $\psi_{ij}$ is an isomorphism}
\]
In addition, we define the submodule $\rad^2_R(M,N)\subset\Hom_R(M,N)$. The submodule $\rad^2_R(M,N)$ consists of 
morphisms $\psi:M\rightarrow N$ such that $\psi$ decomposes as $\psi=g\circ f$,
{\footnotesize\[\xymatrix@C=10pt@R=10pt{
   M\ar[rr]^{\psi}\ar[dr]_f&&N\\
   &Z\ar[ur]_g& \\
}\]}
where $Z$ is an MCM $R$-module,\;$f\in\rad_R(M,Z)$,\;$g\in\rad_R(Z,N)$.
We call a morphism $\psi:M\rightarrow N$ irreducible if $\psi\in\rad_R(M,N)\setminus\rad^2_R(M,N)$.
Set 
\[
 \Irr_R(M,N)\coloneqq\rad_R(M,N)\big/\rad^2_R(M,N),
\]
then $\Irr_R(M,N)$ is a vector space over $k$.
\end{definition}

By using these notions, we define the Auslander-Reiten quiver.

\begin{definition}[Auslander-Reiten quiver]
The Auslander-Reiten $($= AR$)$ quiver of $R$ is an oriented graph whose vertices are indecomposable MCM $R$-modules $R,M_1,\cdots,M_{n-1}$ with 
$\dim_k\Irr_R(M_s,M_t)$ arrows from $M_s$ to $M_t$ \;$(s,t=0,1,\cdots,n-1)$. 
\end{definition}

We will give an example of an AR quiver in the next subsection. 

%%%%%%%%%%%%%%%%%%%%%%%%%%%%%%%%%%%%%%%%%
\subsection{The case of cyclic quotient surface singularities} 
\label{subsec_cyclic}

Since one of the aims of this paper is to determine the dual $F$-signature of special CM modules over cyclic quotient surface singularities, 
we restate results in subsection~\ref{subsec_genFsig} and \ref{subsec_AR} for the cyclic case. 
Thus, we suppose that $G$ is a cyclic group as follows.
\[
G\coloneqq \langle\;\sigma=
    \begin{pmatrix} \zeta_n&0 \\
                    0&\zeta_n^a 
    \end{pmatrix}           \;\rangle,
\]
where $\zeta_n$ is a primitive $n$-th root of unity, $1\le a\le n-1$,\;and $\mathrm{gcd}(a,n)=1$.
We denote the cyclic group $G$ as above by $\frac{1}{n}(1,a)$. 
Let $S\coloneqq k[[x,y]]$ be a power series ring and we assume that $n$ is coprime to $p=\mathrm{char}\;k$.
We denote the invariant subring of $S$ under the action of $G$ by $R\coloneqq S^G$. 
Since $G$ is an abelian group, any irreducible representations of $G$ are described by
\[
 V_t:\sigma\mapsto\zeta_n^{-t} \quad(t=0,1,\cdots,n-1).
\]
Then we set, 
\[
 M_t\coloneqq (S\otimes_kV_t)^G=\Bigl\{\sum_{i,j}a_{ij}x^iy^j\in S\;|\;a_{ij}\in k,\;i+ja\equiv t\;(\mathrm{mod}\;n)\Bigl\},\;\;(t=0,1,\cdots,n-1)
\]
%Then, these $M_t$ are finitely many MCM modules over $R$ and $\rank M_t=1$. 
These give all indecomposable MCM modules over $R$, and each has rank one. 

From Corollary~\ref{gene-Fsig-cor}, \;$s(M_t,M_u)=1/n \;(u=0,1,\cdots,n-1)$. 
Thus, when we discuss the asymptotic behavior of ${}^eM_t$ on the order of $p^{2e}$, we may consider as 
\begin{equation}
\label{on_order_p}
{}^eM_t\approx (R\oplus M_1\oplus\cdots\oplus M_{n-1})^{\oplus\frac{p^{2e}}{n}} .
\end{equation}

\bigskip

Also, the AR sequence ending in $M_t\;(t\neq 0)$ is 
\begin{equation}
\label{AR_M}
0\longrightarrow M_{t-a-1}\longrightarrow M_{t-1}\oplus M_{t-a}\longrightarrow M_t\longrightarrow 0.
\end{equation}
In the case of $t=0$, the fundamental sequence of $R$ is
\begin{equation}
\label{fund_R}
0\longrightarrow\omega_R\longrightarrow M_{-1}\oplus M_{-a}\longrightarrow R\longrightarrow k\longrightarrow 0.
\end{equation}

Thus, we have $\tau(M_t)=M_{t-a-1}$ and $E_{M_t}=M_{t-1}\oplus M_{t-a}$ for $t=0,1,\cdots,n-1$. 
It is known that $\dim_k\Irr_R(M_s,M_t)$ is equal to the multiplicity of $M_s$ in the decomposition of $E_{M_t}$.
Therefore, by $(\ref{AR_M})$ and $(\ref{fund_R})$, there is an arrow from $M_{t-1}$ to $M_t$, and from $M_{t-a}$ to $M_t$ for $t=0,1,\cdots,n-1$.  
Namely, we have $\dim_k\Irr_R(M_{t-1},M_t)=1$ and $\dim_k\Irr_R(M_{t-a},M_t)=1$. 

\begin{remark}
\label{xy_surj}
Since $S\cong R\oplus M_1\oplus\cdots\oplus M_{n-1}$, each MCM $R$-modules $M_t$ is an $R$-submodule of $S$, 
and we can take a morphism $\cdot x$\;(resp.\,$\cdot y$) as a basis of $1$-dimensional vector space $\Irr_R(M_{t-1},M_t)$\;(resp. $\Irr_R(M_{t-a},M_t)$).
\[
M_{t-1}=\big\{\;f\in S\;|\;\sigma\cdot f=\zeta_n^{t-1}\;f\;\big\}\overset{x}{\longrightarrow}
M_t=\big\{\;f\in S\;|\;\sigma\cdot f=\zeta_n^t\;f\;\big\} 
\]
\[
M_{t-a}=\big\{\;f\in S\;|\;\sigma\cdot f=\zeta_n^{t-a}\;f\;\big\}\overset{y}{\longrightarrow}
M_t=\big\{\;f\in S\;|\;\sigma\cdot f=\zeta_n^t\;f\;\big\} 
\]
\end{remark}

\begin{example}
Let $G=\frac{1}{7}(1,3)$ be a cyclic group of order $7$.
Irreducible representations of $G$ are   
\[
V_t:\sigma\mapsto\zeta_7^{-t} \quad(t=0,\cdots,6),
\]
where $\zeta_7$ is a primitive $7$-th root of unity.
Then the AR quiver of $R$ is described as follows. 
For simplicity, we only describe subscripts as vertices, and all common numbers are identified. 
{\small
\[\xymatrix@C=15pt@R=12pt{
  0\ar[r]&3\ar[r]&6\ar[r]&\cdots\ar[r]&1\ar[r]&4\ar[r]&0 \\
  6\ar[u]\ar[r]&2\ar[u]\ar[r]&5\ar[u]\ar[r]&\cdots\ar[r]&0\ar[u]\ar[r]&3\ar[u]\ar[r]&6\ar[u] \\
  5\ar[u]\ar[r]&1\ar[u]\ar[r]&4\ar[u]\ar[r]&\cdots\ar[r]&6\ar[u]\ar[r]&2\ar[u]\ar[r]&5\ar[u] \\
  \vdots\ar[u]&\vdots\ar[u]&\vdots\ar[u]&\reflectbox{$\ddots$}&\vdots\ar[u]&\vdots\ar[u]&\vdots\ar[u] \\
  2\ar[u]\ar[r]&5\ar[u]\ar[r]&1\ar[u]\ar[r]&\cdots\ar[r]&3\ar[u]\ar[r]&6\ar[u]\ar[r]&2\ar[u] \\
  1\ar[u]\ar[r]&4\ar[u]\ar[r]&0\ar[u]\ar[r]&\cdots\ar[r]&2\ar[u]\ar[r]&5\ar[u]\ar[r]&1\ar[u] \\
  0\ar[u]\ar[r]&3\ar[u]\ar[r]&6\ar[u]\ar[r]&\cdots\ar[r]&1\ar[u]\ar[r]&4\ar[u]\ar[r]&0\ar[u]  
}\]}
\end{example}

\begin{remark}
For each diagram $\begin{array}{c}\tiny{\xymatrix@C=8pt@R=8pt{a\ar[r]&b\\c\ar[u]\ar[r]&d\ar[u]}}\end{array}$, if $b\neq 0$ then 
$0\rightarrow M_c\rightarrow M_a\oplus M_d\rightarrow M_b\rightarrow 0$ 
is the AR sequence ending in $M_b$, and any diagram commutes  
$\begin{array}{c}\tiny{\xymatrix@C=8pt@R=8pt{a\ar[r]^y\ar@{}[dr]|\circlearrowleft&b\\c\ar[u]^x\ar[r]_y&d\ar[u]_x}}\end{array}$
by Remark~\ref{xy_surj}.
\end{remark}

%%%%%%%%%%%%%%%%%%%%%%%%%%%%%%%%%%%%%%%%%%%%%%%%%%%%%%%%%%%%%%%%%%%%%%%%%%%%%%%%%%%%%%%%%%%%%%%%%%%%%%%%%%%%
\section{Dual $F$-signature of Special CM modules}
\label{section_DFsig}
In this section, we introduce the notion of special CM modules and determine the dual $F$-signature of them. 
Firstly, we recall the definition of special CM modules over an invariant subring $R$, and the properties of them.

\begin{definition}[\cite{Wun2}]
For an MCM $R$-module $M$, we call $M$ special if $(M\otimes_R\omega_R)\big/\tor$ is also an MCM $R$-module.
\end{definition}

In other words, let $\varphi$ be the natural morphism $M\otimes_R\omega_R\rightarrow(M\otimes_R\omega_R)^{**}$,  
then $M\otimes_R\omega_R\big/\Ker\varphi$ is also an MCM $R$-module if and only if $M$ is a special CM $R$-module. 
In that case, we have the following (cf. \cite[Lemma~9]{Rie}), 
\[
M\otimes_R\omega_R\big/\Ker\varphi\cong\tau(M)\cong(M\otimes_R\omega_R)^{**}.\]
Therefore, $M$ is a special CM $R$-module if and only if $\varphi$ is a surjection. 
Furthermore, there are several characterizations of special CM modules as follows (see \cite[Theorem~2.7 and 3.6]{IW}). 

\begin{proposition} 
\label{char_special}
Suppose that $M$ is an MCM $R$-module. Then the following are equivalent.  
\begin{itemize}
 \item[(1)] $M$ is a special CM module,
 \item[(2)] $\Ext^1_R(M,R)=0$, 
 \item[(3)] $(\Omega M)^*\cong M$ where $\Omega M$ is the syzygy of $M$. 
\end{itemize}
\end{proposition}

Suppose $M$ is a special CM $R$-module, then we have the following exact sequence by the condition (3). 
Here, $\mu_R(M)$ is the number of minimal generators of $M$. 
\[
0\rightarrow M^*\cong\Omega M\rightarrow R^{\oplus\mu_R(M)}\rightarrow M\rightarrow 0. 
\]
Thus, we have $\mu_R(M)=2\rank_RM$. The converse is true if $\rank_RM=1$ (cf. \cite[Theorem\,2.1]{Wun2}). 
If $\rank_RM>1$, the converse is no longer true (cf. \cite[Example\,A.5]{Nak} and \cite{IW}).  
Since each MCM module over cyclic quotient surface singularities has rank one, a special CM module is minimally $2$-generated 
(see Theorem~\ref{cyclic_special}).

\bigskip
For a cyclic group $G=\frac{1}{n}(1,a)$, we can describe special CM-modules as follows.

Firstly, we consider the Hirzebruch-Jung continued fraction expansion of $n/a$, 
\[
 \frac{n}{a}=\alpha_1-\cfrac{1}{\alpha_2-\cfrac{1}{\cdots -\frac{1}{\alpha_r}}}\coloneqq[\alpha_1,\alpha_2,\cdots,\alpha_r],
\]
and then we introduce the notion of $i$-series and $j$-series (cf. \cite{Wem}, \cite{Wun1}).

\begin{definition}
For $n/a=[\alpha_1,\alpha_2,\cdots,\alpha_r]$, we define the $i$-series and the $j$-series as follows.
\begin{eqnarray*}
 i_0=n,\;\;i_1=a,\;\;&i_t=\alpha_{t-1}i_{t-1}-i_{t-2}&\;\;(t=2,\cdots,r+1), \\
 j_0=0,\;\;j_1=1,\;\;&j_t=\alpha_{t-1}j_{t-1}-j_{t-2}&\;\;(t=2,\cdots,r+1).
\end{eqnarray*}
\end{definition}

\begin{remark}
\label{rem_i_j}
By the construction method of the $i$-series and the $j$-series, it is easy to see 
 \begin{eqnarray*}
  &\cdot& i_t\equiv j_ta\;(\mathrm{mod}\;n),\\
  &\cdot& i_0=n>i_1=a>i_2>\cdots>i_r=1>i_{r+1}=0,\\
  &\cdot& j_0=0<j_1=1<j_2=\alpha_1<\cdots<j_r<j_{r+1}=n.
%  &\cdot& i_t/i_{t+1}=\alpha_{t+1}-\cfrac{1}{\alpha_{t+2}-\cfrac{1}{\alpha_{t+3}-\cdots}}=[\alpha_{t+1},\alpha_{t+2},\cdots,\alpha_r],\\
% &\cdot& j_{t+1}/j_t=\alpha_t-\cfrac{1}{\alpha_{t-1}-\cfrac{1}{\alpha_{t-2}-\cdots}}=[\alpha_t,\alpha_{t-1},\cdots,\alpha_1] .
 \end{eqnarray*}
\end{remark}

By using the $i$-series and the $j$-series, we can characterize special CM $R$-modules.

\begin{theorem}[\cite{Wun1}]
\label{cyclic_special}
For a cyclic group $G=\frac{1}{n}(1,a)$ with $n/a=[\alpha_1,\alpha_2,\cdots,\alpha_r]$, special CM $R$-modules are
$M_{i_t}\;(t=0,1,\cdots,r)$. Moreover, minimal generators of $M_{i_t}$ are $x^{i_t}$ and $y^{j_t}$ for $t=1,\cdots,r$. 
\end{theorem}

\begin{example}
\label{ex7/3}
Let $G=\frac{1}{7}(1,3)$ be a cyclic group of order $7$. 
The Hirzebruch-Jung continued fraction expansion of $7/3$ is 
\[
 \frac{7}{3}=3-\cfrac{1}{2-1/2}=[3,2,2],
\]
and the $i$-series and the $j$-series are described as follows.
\[
\begin{array}{ccccc}
 i_0=7,&i_1=3,&i_2=2,&i_3=1,&i_4=0, \\
 j_0=0,&j_1=1,&j_2=3,&j_3=5,&j_4=7.
\end{array}
\]
Thus, special CM modules are $R,M_1,M_2,M_3$ and these are described explicitly 
\[
\begin{array}{c}
  R=k[[x^7,x^4y,xy^2,y^7]] \\
  M_1=Rx+Ry^5 \\
  M_2=Rx^2+Ry^3 \\
  M_3=Rx^3+Ry .
\end{array}\]
\end{example}

%By using the AR quiver, we can reinterpret these results as follows.
We now show, using AR theory, how to investigate possible surjections ${}^eM_2\twoheadrightarrow M_2^{\oplus b_e}$. 

{\small
\begin{figure}[!h]
\begin{tabular}{c}
\begin{minipage}{0.5\hsize}
\[\xymatrix@C=10pt@R=10pt{
  \ar@{.}[r]&4\ar@{>}[r]&0\ar@{=>}[r]^y&3\ar@{=>}[r]^y&6\ar@{=>}[r]^y&*++[Fo]{2} \\
  \ar@{.}[r]&3\ar@{>}[u]\ar@{>}[r]&6\ar@{>}[u]\ar@{>}[r]&2\ar@{>}[u]\ar@{>}[r]&5\ar@{>}[u]\ar@{>}[r]&1\ar@{=>}[u]_x \\
  \ar@{.}[r]&2\ar@{>}[u]\ar@{>}[r]&5\ar@{>}[u]\ar@{>}[r]&1\ar@{>}[u]\ar@{>}[r]&4\ar@{>}[u]\ar@{>}[r]&0\ar@{=>}[u]_x \\
  \ar@{.}[r]&1\ar@{>}[u]\ar@{>}[r]&4\ar@{>}[u]\ar@{>}[r]&0\ar@{>}[u]\ar@{>}[r]&3\ar@{>}[u]\ar@{>}[r]&6\ar@{>}[u] \\
 &\ar@{.}[u]&\ar@{.}[u]&\ar@{.}[u]&\ar@{.}[u]&\ar@{.}[u] }\]
\caption{}
\label{explain_surj_1} 
\end{minipage} 

\begin{minipage}{0.5\hsize}
\[\xymatrix@C=10pt@R=10pt{
  \ar@{.}[r]&\bullet\ar@{>}[r]&0\ar[r]^y&3\ar[r]^y&6\ar[r]^y&*++[Fo]{2} \\
  \ar@{.}[r]&\bullet\ar@{>}[u]\ar@{>}[r]&\bullet\ar@{>}[u]\ar@{>}[r]&2\ar@{>}[u]\ar@{>}[r]&5\ar@{>}[u]\ar@{>}[r]&1\ar[u]_x \\
  \ar@{.}[r]&\bullet\ar@{>}[u]\ar@{>}[r]&\bullet\ar@{>}[u]\ar@{>}[r]&\bullet\ar@{>}[u]\ar@{>}[r]&\bullet\ar@{>}[u]\ar@{>}[r]&0\ar[u]_x \\
  \ar@{.}[r]&\bullet\ar@{>}[u]\ar@{>}[r]&\bullet\ar@{>}[u]\ar@{>}[r]&\bullet\ar@{>}[u]\ar@{>}[r]&\bullet\ar@{>}[u]\ar@{>}[r]&\bullet\ar@{>}[u] \\
 &\ar@{.}[u]&\ar@{.}[u]&\ar@{.}[u]&\ar@{.}[u]&\ar@{.}[u] }\]
\caption{}
\label{explain_surj_2}
\end{minipage} 
\end{tabular}
\end{figure} }

We take the MCM $R$-module $M_2$ as an example.
From the AR quiver around the vertex $\textcircled{\footnotesize2}$, we can see that there are several morphisms ending in $\textcircled{\footnotesize2}$ 
and obtain minimal generators $x^2$ and $y^3$ by following the morphisms described by double arrows in Figure~\ref{explain_surj_1}. 

Since each diagrams $\begin{array}{c}\tiny{\xymatrix@C=7pt@R=7pt{a\ar[r]&b\\c\ar[u]\ar[r]&d\ar[u]}}\end{array}$ are commute, morphisms from vertices 
which are denoted by $\bullet$ in Figure~\ref{explain_surj_2} to $\textcircled{\footnotesize2}$ go through $``0"(\text{that is}, ``R")$.
Thus, the image of each morphism $\bullet\rightarrow\textcircled{\footnotesize2}$ is in $\fkm M_2$ where $\fkm$ is the maximal ideal of $R$. 
By Nakayama's lemma, such a morphism doesn't contribute to construct a surjection. Thus we may ignore them. 
Also, there are morphisms from vertices which are denoted by $\bigstar$ in Figure~\ref{explain_surj_3} to $\textcircled{\footnotesize2}$. 
Minimal generators of $M_\bigstar$ are generated by morphisms from $0$ $($which are located outside of dotted area in Figure~\ref{explain_surj_3}$)$ to $\bigstar$. 
Considering the composition of such a morphism and $\bigstar\rightarrow\textcircled{\footnotesize2}$
\[
 R\rightarrow M_\bigstar\rightarrow M_2\;(1\mapsto\delta\mapsto x^{m_1}y^{m_2}\delta),
\]
where $\delta$ is a minimal generator of $M_\bigstar$ and $m_1\ge1,\;m_2\ge1$. 
Then it is easy to see that the image of the morphism $\bigstar\rightarrow\textcircled{\footnotesize2}$ is in $\fkm M_2$. Thus we may ignore them. 

{\small
\begin{figure}[!h]
\begin{tabular}{c}
\begin{minipage}{0.5\hsize}
\[\xymatrix@C=10pt@R=10pt{
  &&0\ar[r]^y\ar@{.}[dd]&3\ar[r]^y&6\ar[r]^y&*++[Fo]{2} \\
  0\ar@{-->}[rrr]&&&\bigstar\ar@{>}[u]_x\ar@{>}[r]&\bigstar\ar@{>}[u]_x\ar@{>}[r]&1\ar[u]_x \\
  &&\ar@{.}[rrr]&&&0\ar[u]_x \\
  &&&0\ar@{-->}[uu]&&\\
  &&&&0\ar@{-->}[uuu]&
}\]
\caption{}
\label{explain_surj_3}
\end{minipage} 

\begin{minipage}{0.5\hsize}
\[\xymatrix@C=10pt@R=10pt{
  0\ar[r]^y&3\ar[r]^y&6\ar[r]^y&*++[Fo]{2} \\
  &&&1\ar[u]_x \\
  &&&0\ar[u]_x 
}\]
\caption{}
\label{explain_surj_4}
\end{minipage} 
\end{tabular}
\end{figure} }

Thus, in order to investigate a surjection ${}^eM_2\twoheadrightarrow M_2^{\oplus b_e}$, we need only discuss the MCM $R$-modules located in 
the horizontal direction from $M_2$ to $R$ and the vertical direction from $M_2$ to $R$ $($Figure~\ref{explain_surj_4}$)$.

\medskip

In general, the number of minimal generators of a special CM $R$-module $M_{i_t}$ is two and minimal generators take a form like $x^{i_t},\;y^{j_t}$ 
by Theorem~\ref{cyclic_special}. 
Thus, it is equivalent to there is no $``0"$ in dotted vertices area of Figure~\ref{explain_surj_5}. 
By the above arguments, in order to construct a surjection ${}^eM_{i_t}\twoheadrightarrow M_{i_t}^{\oplus b_e}$, we may only discuss 
horizontal direction arrows from $R$ to $M_{i_t}$ and vertical direction arrows from $R$ to $M_{i_t}$. 
We consider sets of subscripts of vertices $\calF_t=\{0,1,\cdots,i_t-1\}$ and $\calG_t=\{i_t-a,\cdots,i_t-j_ta\equiv 0\}$ as in Figure~\ref{explain_surj_5}. 
It is easy to see that $|\calF_t|=i_t,\; |\calG_t|=j_t$.

{\small
\begin{figure}[!h]
\[\xymatrix@C=17pt@R=15pt{
\ar@{-} `u[r] `[rrrr]^*+{\calG_t}&&&&&& \\
0\equiv i_t-j_ta\ar[r]^(0.45){y}&i_t-(j_t-1)a\ar[r]^(0.65){y}&\cdots\ar[r]^(0.4){y}&i_t-2a\ar[r]^y&i_t-a\ar[r]^(0.55){y}&i_t& \\
\bullet\ar@{>}[u]\ar@{>}[r]&\bullet\ar@{>}[u]\ar@{>}[r]&\cdots\ar[r]&\bullet\ar@{>}[u]\ar@{>}[r]&\bullet\ar@{>}[u]\ar@{>}[r]&i_t-1\ar[u]_x &\\
\vdots\ar@{>}[u]&\vdots\ar@{>}[u]&\reflectbox{$\ddots$}&\vdots\ar@{>}[u]&\vdots\ar@{>}[u]&\vdots\ar[u]_x &\\
\bullet\ar@{>}[u]\ar@{>}[r]&\bullet\ar@{>}[u]\ar@{>}[r]&\cdots\ar[r]&\bullet\ar@{>}[u]\ar@{>}[r]&\bullet\ar@{>}[u]\ar@{>}[r]&1\ar[u]_(0.45){x}&\\
\bullet\ar@{>}[u]\ar@{>}[r]&\bullet\ar@{>}[u]\ar@{>}[r]&\cdots\ar[r]&\bullet\ar@{>}[u]\ar@{>}[r]&\bullet\ar@{>}[u]\ar@{>}[r]&0\ar[u]_x&\ar@{-} `r[uu] `[uuu]_*+{\calF_t}
}\]
\caption{}
\label{explain_surj_5}
\end{figure} }

To determine the dual $F$-signature of special CM $R$-modules, we prepare some notations and lemmas. 

For the $i$-series $(i_1,\cdots,i_r)$ associated with $\frac{1}{n}(1,a)$ and any $\beta\in\mathbb{Z}_{\ge 0}$ with $0\le\beta\le n-1$, 
there are unique non-negative integers $d_1,\cdots,d_r\in\mathbb{Z}_{\ge 0}$ such that 
\begin{eqnarray*}
  \beta=d_1i_1+h_1,\quad &h_1\in\mathbb{Z}_{\ge 0},&\quad 0\le h_1<i_1, \\
  h_t=d_{t+1}i_{t+1}+h_{t+1},\quad &h_{t+1}\in\mathbb{Z}_{\ge 0},&\quad 0\le h_{t+1}<i_{t+1},\quad(t=1,\cdots,r-1), \\
  &h_r=0.&
\end{eqnarray*}

Thus, we can describe $\beta$ as follows, 
\[
  \beta=d_1i_1+d_2i_2+\cdots+d_ri_r.
\]
For such $\beta$, there is the unique integer $\widetilde{\beta}\in\mathbb{Z}_{\ge 0}$ such that $a\widetilde{\beta}\equiv\beta\;(\mathrm{mod}\;n),
\;0\le\widetilde{\beta}\le n-1$. 

\begin{lemma}[\cite{Wun1}]
\label{d-j-series}
Let $\widetilde{\beta}$ be the same as above. Then $\widetilde{\beta}$ is described as 
\[
 \widetilde{\beta}=d_1j_1+d_2j_2+\cdots+d_rj_r,
\]
where $(j_1,\cdots,j_r)$ is the $j$-series associated with $\frac{1}{n}(1,a)$. 
\end{lemma}

\begin{lemma}
\label{F_intersection_G}
Let the notation be the same as above, then $\calF_t\cap\calG_t=\{0\}$ as a set of subscripts of vertices.
\end{lemma}

\begin{proof}
It is trivial that $0\in\calF_t\cap\calG_t$ by the definition of $\calF_t$ and $\calG_t$. Thus, it suffices to show there is no 
pair $(m_1,m_2)\in\mathbb{Z}_{> 0}^2$ such that $m_1\equiv m_2a \;(\mathrm{mod}\;n)$, where $1\le m_1\le i_t-1$ and $1\le m_2\le j_t-1$. 
Assume that there exists such a pair $(m_1,m_2)$. Then there are non-negative integers $d_1,\cdots,d_r$ such that $m_1=d_1i_1+d_2i_2+\cdots+d_ri_r$. 
Since $1\le m_1\le i_t-1$ and $i_t>i_{t+1}$ (cf. Remark~\ref{rem_i_j}), $d_1=\cdots=d_t=0$ and there exists $\lambda$ such that $t+1\le\lambda\le r$ 
and $d_\lambda\neq 0$. 
From Lemma~\ref{d-j-series}, we obtain $m_2=d_1j_1+d_2j_2+\cdots+d_rj_r$. Thus, 
\[
 m_2=d_{t+1}j_{t+1}+\cdots+d_rj_r\ge j_\lambda>j_t.
\]
This contradicts $m_2\le j_t-1$.
\end{proof}

We are now ready to state the main theorem. 

\begin{theorem}
\label{main}
Let the notation be the same as above, then for any non-free special CM $R$-module $M_{i_t}$ one has 

\begin{equation*}
s(M_{i_t})=\begin{cases} \displaystyle\frac{\operatorname{min}(i_t,j_t)+1}{n}&\text{(if \;$i_t\neq j_t$)} \\
                         &\\
                         \displaystyle\frac{2i_t+1}{2n}&\text{(if \;$i_t=j_t$)}. 
           \end{cases}
\end{equation*}
\end{theorem}

\medskip

\begin{proof}

In order to determine the value of the dual $F$-signature of $M_{i_t}$, 
we have to find the maximum number $b_e$ such that there is a surjection ${}^eM_{i_t}\twoheadrightarrow M_{i_t}^{\oplus b_e}$. 
Note that we may consider ${}^eM_{i_t}$ as 
\[
{}^eM_{i_t}\approx (R\oplus M_1\oplus\cdots\oplus M_{n-1})^{\oplus\frac{p^{2e}}{n}}
\]
by (\ref{on_order_p}), hence we may assume the number of each indecomposable MCM module in ${}^eM_{i_t}$ is the same on the order of $p^{2e}$. 
%Thus, we will discuss a surjection from $R\oplus M_1\oplus\cdots\oplus M_{n-1}$ to a finite direct sum of some copies of $M_{i_t}$. 
Let $\calF_t$, $\calG_t$ be the sets of vertices as in Figure~\ref{explain_surj_5}. 
By the above observations, MCM modules which contribute to construct a surjection are $M_{i_t}$ itself 
and modules corresponding to elements in $\calF_t$ or $\calG_t$. 
Since an indecomposable MCM module which is not isomorphic to $R$ and $M_{i_t}$ could construct at most one generator of $M_{i_t}$, 
we should first combine MCM modules corresponding to elements in $\calF_t\setminus\{0\}$ with those in $\calG_t\setminus\{0\}$ for constructing surjections efficiently, 
and then we should use $R$ and $M_{i_t}$. 
Therefore, in what follows, we will find disjoint sets of summands of ${}^eM_{i_t}$ which surject onto $M_{i_t}$ as much as possible along this strategy. 

Firstly, we show the case of $i_t>j_t$. Thus, $|\calF_t|>|\calG_t|$. 
We choose elements $f_1$ and $g_1$ from $\calF_t\setminus\{0\}$ and $\calG_t\setminus\{0\}$ respectively, 
and consider corresponding indecomposable MCM $R$-modules $M_{f_1}$ and $M_{g_1}$. Here, we remark that $f_1\neq g_1$ by Lemma~\ref{F_intersection_G}. 
Then we can construct a surjection $M_{f_1}\oplus M_{g_1}\twoheadrightarrow M_{i_t}$. 
\[
\footnotesize{\xymatrix@C=12pt@R=12pt{
0\ar[r]&\cdots\ar[r]&g_1\ar[r]&\cdots\ar[r]&i_t \\
&&&&\vdots\ar[u] \\
&&&&f_1\ar[u] \\
&&&&\vdots\ar[u] \\
&&&&0\ar[u] 
}}\]
Then we consider the sets $\calF_t\setminus\{0,f_1\}$ and $\calG_t\setminus\{0,g_1\}$. 
Similarly, we choose elements $f_2$ and $g_2$ from the sets $\calF_t\setminus\{0,f_1\}$ and $\calG_t\setminus\{0,g_1\}$ respectively, and construct a surjection $M_{f_2}\oplus M_{g_2}\twoheadrightarrow M_{i_t}$. 
By repeating the same process, we finally arrive at $\calG_t\setminus\{0,g_1, \cdots, g_{j_t-1}\}=\emptyset$ and have $j_t-1$ surjections. 
Since we still don't use $0\in\calG_t$ (that is $R$), we construct a surjection by combining $R$ and an indecomposable MCM module corresponding to 
an element $f^\prime\in\calF_t\setminus\{0,f_1,\cdots, f_{j_t-1}\}\neq\emptyset$.  
In addition, there is a trivial surjection $M_{i_t}\twoheadrightarrow M_{i_t}$. 
Thus, through these processes, we could obtain disjoint sets of summands 
\[
\{M_{f_1}, M_{g_1}\}, \cdots, \{M_{f_{j_t-1}}, M_{g_{j_t-1}}\}, \{M_{f^\prime}, R\}, \{M_{i_t}\}
\] which surject onto $M_{i_t}$. 
Thus, we have a surjection 
\[(R\oplus M_1\oplus\cdots\oplus M_{n-1})^{\oplus\frac{p^{2e}}{n}}\twoheadrightarrow M_{i_t}^{\oplus \frac{p^{2e}}{n}(j_t+1)}.\] 
Therefore the dual $F$-signature of $M_{i_t}$ is $s(M_{i_t})=\displaystyle\frac{j_t+1}{n}$. 

Similarly, we obtain $s(M_{i_t})=\displaystyle\frac{i_t+1}{n}$ for the case of $i_t<j_t$.

For the case of $i_t=j_t$, we repeat the same process until we have $\calF_t\setminus\{0,f_1, \cdots, f_{i_t-1}\}=\emptyset$ and $\calG_t\setminus\{0,g_1, \cdots, g_{j_t-1}\}=\emptyset$, and then we have a surjection 
\[
(M_1\oplus\cdots\oplus M_{i_t-1} \oplus M_{i_t+1}\oplus\cdots\oplus M_{n-1})^{\oplus\frac{p^{2e}}{n}}\twoheadrightarrow M_{i_t}^{\oplus \frac{p^{2e}}{n}(i_t-1)}. 
\]
In addition, there is a trivial surjection $M_{i_t}\twoheadrightarrow M_{i_t}$. 
For now, we don't use $R$, and by using two free summands we also construct the surjection: 
\[
R\oplus R\overset{(x^{i_t}\,\,y^{j_t})}{\twoheadrightarrow} M_{i_t}. 
\]
Thus, the dual $F$-signature of $M_{i_t}$ is 
\[
s(M_{i_t})=\frac{i_t-1}{n}+\frac{1}{n}+\frac{1}{2n}=\frac{2i_t+1}{2n}.
\]
\end{proof}

%   \medskip
   
\begin{example}
Let the notation be as in Example~\ref{ex7/3}. Then, the dual $F$-signature of special CM modules are
\[
 s(M_1)=\frac{2}{7},\;s(M_2)=\frac{3}{7},\;s(M_3)=\frac{2}{7}.
\]
\end{example}

\medskip

Next, we give an example for the case $i_t=j_t$.
\begin{example}
\label{ex8/5}
Let $G=\frac{1}{8}(1,5)$ be a cyclic group of order $8$. 
The Hirzebruch-Jung continued fraction expansion of $8/5$ is 
\[
 \frac{8}{5}=2-\cfrac{1}{3-1/2}=[2,3,2],
\]
and the $i$-series and the $j$-series are described as follows.
\[
\begin{array}{ccccc}
 i_0=8,&i_1=5,&i_2=2,&i_3=1,&i_4=0, \\
 j_0=0,&j_1=1,&j_2=2,&j_3=5,&j_4=8.
\end{array}
\]
Thus, special CM modules are $R,M_1,M_2,M_5$. In this case, we have $i_2=j_2$ and there are surjections as follows.
\[
\begin{array}{cc}
 \begin{array}{c}
 \xymatrix@C=10pt@R=10pt{
  0\ar[r]^y&5\ar[r]^y&2 &\\
  &&1\ar[u]_x &\\
  &&0\ar[u]_x &}
 \end{array}&
 \begin{array}{lcl}
  M_2&\twoheadrightarrow&M_2\\
  M_1\oplus M_5&\twoheadrightarrow&M_2\\
  R\oplus R&\twoheadrightarrow&M_2.
 \end{array}
\end{array}
\]
   
   \medskip
Thus, the dual $F$-signature of $M_2$ is 
\[
s(M_2)=\frac{1}{8}+\frac{1}{8}+\frac{1}{16}=\frac{5}{16}.
\]
\end{example}

\begin{example}
\label{exAtype}
Let $G=\frac{1}{n}(1,n-1)\subset\SL(2,k)$ be a cyclic group of order $n$, that is, Dynkin type $A_{n-1}$.  
The Hirzebruch-Jung continued fraction expansion of $n/(n-1)$ is 
\[
 \frac{n}{n-1}=2-\cfrac{1}{2-\cfrac{1}{\cdots -1/2}}=[\underbrace{2,2,\cdots,2}_{n-1}],
\]
and the $i$-series and the $j$-series are described as follows, 
\[
\begin{array}{llllll}
 i_0=n,&i_1=n-1,&i_2=n-2,&\cdots,&i_{n-1}=1,&i_n=0, \\
 j_0=0,&j_1=1,&j_2=2,&\cdots,&j_{n-1}=n-1,&j_n=n.
\end{array}
\]
Namely, $i_t=n-t,\;j_t=t\;(t=1,2,\cdots,n-1)$.
As we mentioned in Remark~\ref{remSL2}, all irreducible representations of $G=\frac{1}{n}(1,n-1)\subset\SL(2,k)$ are special. 
Thus, any $M_t$ is a special CM module and the dual $F$-signature of $M_t$ is obtained by Theorem~\ref{main}.
\begin{equation*}
s(M_{i_t})=\begin{cases} \displaystyle\frac{1}{n}+\frac{j_t}{n}=\frac{t+1}{n}&\text{(if \;$t<\frac{n}{2}$)} \\
                         &\\
                         \displaystyle\frac{1}{n}+\frac{t-1}{n}+\frac{1}{2n}=\frac{2t+1}{2n}&\text{(if \;$t=\frac{n}{2}$)} \\
                         &\\
                         \displaystyle\frac{1}{n}+\frac{i_t}{n}=\frac{n-t+1}{n}&\text{(if \;$t>\frac{n}{2}$)}. 
           \end{cases}
\end{equation*}
For other Dynkin types (i.e. $D_n,E_6,E_7,E_8$), see \cite{Nak}.
\end{example}

%%%%%%%%%%%%%%%%%%%%%%%%%%%%%%%%%%%%%%%%%%%%%%%%%%%%%%%%%%%%%%%%%%%%%%%%%%%%%%%%%%%%%%%%%%%%%%%%%%%%%%%%%%%%
\section{Comparing with the Auslander-Reiten translation}
\label{compare_AR_translation}

In this section, we compare the dual $F$-signature of a special CM module with its AR translation. 
It will give us a characterization of Gorensteiness (see Theorem~\ref{comp_AR_tr}). 
As we mentioned in Section~\ref{intro}, it is an analogue of Theorem~\ref{charac_sing}\,$(4), (5)$.

The statements appearing in this section are valid for any 
quotient surface singularities. Therefore, we suppose that $G$ is a finite subgroup of $\GL(2,k)$ which contains no pseudo-reflections except the identity, 
and $S\coloneqq k[[x,y]]$ be a power series ring. 
We assume that the order of $G$ is coprime to $p=\mathrm{char}\;k$. We denote the invariant subring of $S$ under the action of $G$ by $R\coloneqq S^G$.
Let $V_0=k,V_1,\cdots,V_n$ be the complete set of irreducible representations of $G$ 
and set indecomposable MCM $R$-modules $M_t\coloneqq (S\otimes_kV_t)^G\;\;(t=0,1,\cdots,n)$.

\begin{lemma}
\label{asymptotic}
Let $M_t$ be an MCM $R$-module as above. Then we have
\begin{equation}
\label{asymptotic_1}
{}^eM_t\approx (R^{\oplus d_{0,t}}\oplus M_1^{\oplus d_{1,t}}\oplus\cdots\oplus M_n^{\oplus d_{n,t}})^{\oplus\frac{p^{2e}}{n}}\approx {}^e\tau(M_t)
\end{equation}
on the order of $p^{2e}$ ($e\gg0$), where $d_{i,t}=(\rank_RM_t)\cdot(\rank_RM_i)$ and $\tau$ stands for the AR translation.
Furthermore, we have
\[
R^{\oplus d_{0,t}}\oplus M_1^{\oplus d_{1,t}}\oplus\cdots\oplus M_n^{\oplus d_{n,t}}\cong 
\tau(R)^{\oplus d_{0,t}}\oplus \tau(M_1)^{\oplus d_{1,t}}\oplus\cdots\oplus \tau(M_n)^{\oplus d_{n,t}}.
\]

\begin{proof}
From Corollary~\ref{gene-Fsig-cor}, we may write 
\[
{}^eM_t\approx (R^{\oplus d_{0,t}}\oplus M_1^{\oplus d_{1,t}}\oplus\cdots\oplus M_n^{\oplus d_{n,t}})^{\oplus\frac{p^{2e}}{n}},
\]
\[
{}^e\tau(M_t)\approx (R^{\oplus d^\prime_{0,t}}\oplus M_1^{\oplus d^\prime_{1,t}}\oplus\cdots\oplus M_n^{\oplus d^\prime_{n,t}})^{\oplus\frac{p^{2e}}{n}},
\]
where $d^\prime_{i,t}=(\rank_R\tau(M_t))\cdot(\rank_RM_i)$.
Since $\rank_RM_t=\rank_R\tau(M_t)$, it follows that $d_{i,t}=d^\prime_{i,t} \;(i=0,1,\cdots,n)$. This implies (\ref{asymptotic_1}).

Since the AR translation $\tau$ gives a bijection from the set of finitely many indecomposable MCM $R$-modules to itself, 
we set $\tau(M_i)=M_{\sigma(i)} \;(i=0,1,\cdots,n)$ where $\sigma$ is an element of symmetric group $\fkS_{n+1}$.
Then we have 
\[
R^{\oplus d_{0,t}}\oplus M_1^{\oplus d_{1,t}}\oplus\cdots\oplus M_n^{\oplus d_{n,t}}=
M_{\sigma(0)}^{\oplus d_{\sigma(0),t}}\oplus M_{\sigma(1)}^{\oplus d_{\sigma(1),t}}\oplus\cdots\oplus M_{\sigma(n)}^{\oplus d_{\sigma(n),t}}, 
\]
and 
\begin{eqnarray*}
d_{\sigma(i),t}&=&(\rank_RM_t)\cdot(\rank_RM_{\sigma(i)})=(\rank_RM_t)\cdot(\rank_R\tau(M_i)) \\
 &=&(\rank_RM_t)\cdot(\rank_RM_i)=d_{i,t}.
\end{eqnarray*}
Thus, 
\[
M_{\sigma(0)}^{\oplus d_{\sigma(0),t}}\oplus M_{\sigma(1)}^{\oplus d_{\sigma(1),t}}\oplus\cdots\oplus M_{\sigma(n)}^{\oplus d_{\sigma(n),t}}= 
\tau(R)^{\oplus d_{0,t}}\oplus \tau(M_1)^{\oplus d_{1,t}}\oplus\cdots\oplus \tau(M_n)^{\oplus d_{n,t}}.
\]
\end{proof}
\end{lemma}

\begin{theorem}
\label{comp_AR_tr}
For any indecomposable special CM $R$-module $M_t$, we have 
\[
s(M_t)\le s(\tau(M_t)).
\]
Moreover, if $s(M_t)= s(\tau(M_t))$ for an indecomposable special CM $R$-module $M_t$, then $R$ is Gorenstein. 
(Note that if $R$ is Gorenstein, then $s(M)=s(\tau(M))$ holds for all indecomposable MCM modules.) 
%Moreover, $R$ is Gorenstein if and only if $s(M_t)= s(\tau(M_t))$ for an indecomposable special CM $R$-module $M_t$. 
\end{theorem}

\begin{proof}
From Lemma~\ref{asymptotic}, we may write 
\[
{}^eM_t\approx {}^e\tau(M_t)\approx(R^{\oplus d_{0,t}}\oplus M_1^{\oplus d_{1,t}}\oplus\cdots\oplus M_n^{\oplus d_{n,t}})^{\oplus\frac{p^{2e}}{n}}
\]
when we discuss the asymptotic behavior on the order of $p^{2e}$ where $d_{i,t}=(\rank_RM_t)\cdot(\rank_RM_i)$.
%In the rest of this proof, we discuss on this setting, and for simplicity we identify ${}^eM_t\approx{}^e\tau(M_t)$ with $R^{\oplus d_{0,t}}\oplus M_1^{\oplus d_{1,t}}\oplus\cdots\oplus M_n^{\oplus d_{n,t}}$. 

Let $b_e\coloneqq b_e(M_t)$ be the $e$-th $F$-surjective number of $M_t$, hence there exists a surjection ${}^eM_t\twoheadrightarrow M_t^{\oplus b_e}$.
Since $M_t$ is special, the number of minimal generators of $M_t$ is equal to $u\coloneqq 2\rank_RM_t$. 
Thus, there exists a surjection $R^{\oplus b_eu}\twoheadrightarrow M_t^{\oplus b_e}$ which induces the following commutative diagram. 
\[
\xymatrix@C=15pt@R=15pt{
{}^eM_t\ar@{>>}[r]&M_t^{\oplus b_e}\\
R^{\oplus b_eu}\ar@{.>}[u]\ar@{>>}[ru]\\
}\]
Applying the functor $(-\otimes_R\omega_R)^{**}$ to this commutative diagram, then we obtain the commutative diagram. 
\[
\xymatrix@C=20pt@R=15pt{
{}^e\tau(M_t)\approx({}^eM_t\otimes_R\omega_R)^{**}\ar[r]^(0.65){\psi_2}&\tau(M_t)^{\oplus b_e}\\
\omega_R^{\oplus b_eu}\ar[u]\ar[ru]_{\psi_1}\\
}\]
Note that the morphism $\psi_1$ is surjective because the surjection $R^{\oplus b_eu}\twoheadrightarrow M_t^{\oplus b_e}$ induces 
\[
\xymatrix@C=25pt@R=20pt{
\omega_R^{\oplus b_eu}\ar@{>>}[r]\ar[d]_{\cong}&(M_t\otimes_R\omega_R)^{\oplus b_e}\ar@{>>}[d]^{\oplus\varphi}\\
((\omega_R)^{**})^{\oplus b_eu}\ar[r]_(0.45){\psi_1}&((M_t\otimes_R\omega_R)^{**})^{\oplus b_e}
}\]
and $\varphi:M_t\otimes_R\omega_R\rightarrow(M_t\otimes_R\omega_R)^{**}$ is surjective. This implies $\psi_2$ is also surjective, and we have $s(M_t)\le s(\tau(M_t))$. 

   \medskip

If $R$ is Gorenstein, then $M\cong\tau(M)$ for all indecomposable MCM modules. Thus $s(M)=s(\tau(M))$ holds. 
In the rest, we assume that $R$ is not Gorenstein, and hence $R\not\cong\omega_R$. 
Therefore we also have $M\not\cong\tau(M)$ for all indecomposable MCM modules. 
For any indecomposable special CM module $M_t$, we have the following surjection by the same way as above
\[
\omega_R^{\oplus b_eu}\longrightarrow{}^e\tau(M_t)\approx \left(R^{\oplus d_{0,t}}\oplus\bigoplus_{i=1}^n M_i^{\oplus d_{i,t}}\right)^{\oplus\frac{p^{2e}}{n}}
\overset{\psi_2}{\twoheadlongrightarrow}\tau(M_t)^{\oplus b_e}. 
\]
In this surjection, morphisms which go through $R$ don't contribute to construct a surjection by Nakayama's lemma. 
Thus, in addition to a surjection $\omega_R^{\oplus b_eu}\twoheadlongrightarrow\tau(M_t)^{\oplus b_e}$, 
%\[
%\left(\bigoplus_{i=0}^n M_i^{\oplus d_{i,t}}\right)^{\oplus\frac{p^{2e}}{n}}\Big/ R^{\oplus \frac{p^{2e}}{n}d_{0,t}}\twoheadlongrightarrow\tau(M_t)^{\oplus b_e}
%\]
%is also surjective. 
we also construct a surjection 
\[
R^{\oplus \frac{p^{2e}}{n}d_{0,t}}\twoheadlongrightarrow\tau(M_t)^{\oplus \frac{p^{2e}}{n}\frac{d_{0,t}}{v}},
\]
where $v$ is the number of minimal generators of $\tau(M_t)$. Therefore, we obtain 
\[
 b_e(\tau(M_t))\ge b_e+\frac{d_{0,t}p^{2e}}{vn}
\]
where $b_e(\tau(M_t))$ is the $e$-th $F$-surjective number of $\tau(M_t)$. Thus, 
\[
s(\tau(M_t))\ge s(M_t)+\frac{d_{0,t}}{vn} >s(M_t).
\]
\end{proof}

\subsection*{Acknowledgements}The author is deeply grateful to Professor Mitsuyasu Hashimoto for giving him valuable advice and encouragements. 
He also thanks Akiyoshi Sannai for giving some comments about the dual $F$-signature and thanks Professor Ken-ichi Yoshida for suggesting the comparison between special CM modules and other modules (Section~\ref{compare_AR_translation} is based on his suggestion). 
Finally, the author would also like to thank the referee for his/her careful reading and useful comments. 

The author is supported by Grant-in-Aid for JSPS Fellows $($No. 26-422$)$.

\end{document}